\documentclass[draft]{amsart}
\usepackage{amssymb,latexsym,amsfonts,verbatim,amscd}

\numberwithin{equation}{section}

\theoremstyle{plain}

\newtheorem{theorem}[equation]{Theorem}   

\newtheorem{proposition}[equation]{Proposition}

\theoremstyle{definition}

\newtheorem{remark}[equation]{Remark}

\DeclareMathOperator{\Tor}{Tor} \DeclareMathOperator{\Hom}{Hom}

\DeclareMathOperator{\rank}{rank}

 \DeclareMathOperator{\modulo}{mod}
\DeclareMathOperator{\card}{card}

\def\mapright#1{\smash{\mathop{\longrightarrow}\limits^{#1}}}

\renewcommand\dim{{\rm dim}}

\newcommand{\bj}{\boldsymbol{j}}
\newcommand{\bff}{\boldsymbol{f}}

\newcommand\G{\Gamma}
\newcommand\D{\Delta}

\newcommand{\e}{\epsilon}

\newcommand\I{{\rm I}}
\newcommand\J{{\rm J}}

\renewcommand\L{{\rm L}}
\newcommand\<{{\rm <}}
\newcommand{\bigger}{{\rm >}}

\begin{document}

\title[The Denominator of Poincar\'e series]
{On the Denominator of the Poincar\'e series for monomial quotient
rings}
\author[H. Charalambous]{Hara Charalambous}
\address{Department of Mathematics\\
         University at Albany, SUNY\\
         Albany, NY 12222}
\email{hara@math.albany.edu} \keywords{} \subjclass{}
\date{\today}

\begin{abstract}
Let $S=\Bbbk[x_1,\ldots, x_n]$ be a polynomial ring over a field
$\Bbbk$ and $I$ a monomial ideal of $S$. It is well known that the
Poincar\'e series of $\Bbbk$ over $S/I$ is rational. We describe
the coefficients of the denominator of the series and study the
multigraded homotopy Lie algebra of $S/I$.
\end{abstract}

\maketitle


\tableofcontents

\section{Introduction}\label{S: introd}

Let $S=\Bbbk[x_1,\ldots, x_n]$ be a polynomial ring over a field
$\Bbbk$, $I$ a monomial ideal of $S$ and $R=S/I$. If $\J$ is a
subset of the minimal  monomial generating set $\I$ of $I$, we let
$m_\J$ denote the least common multiple of the monomials in $\J$
or the corresponding monomials in $\Bbbk[y_1,\ldots, y_n]$ as
appropriate. The multigraded Poincar\'e series of $R$, $P_R(y,t)$
is  rational, \cite{Ba82}:
 $$P_R(y,t) = \sum_{i=0}^{\infty}  \sum_{{\bf j}\ge 0}
y^{\bf j} {\dim}_k ({\Tor}_i^R(k,k)_{\bf j})
t^i={{\prod_{i=1}^n(1+ty_i)}\over{Q_R(y,t)}}$$ where $Q_R(y,t)\in
Z[y_1,\ldots, y_n][t]$.

The few facts that are known in general about the multigraded
expansion $Q_R(y,t)$ also date back to \cite{Ba82}: the degree of
$Q_R(y,t)$ in $t$ is bounded above by the total degree of $m_\I$,
and the  monomial coefficients of $t$ divide $m_\I$. In this paper
we  discuss the monomial coefficients of $Q_R(y,t)$.

In the first section  of this paper we show that the monomial
coefficients of $t$ in $Q_R(y,t)$ are  least common multiples of
the monomial generators of $I$. We discuss the Koszul homology of
$R$ when $R$ is Golod. We also present a free resolution of $R$
when $I$ is generic.
 In the second section we  discuss the multigraded acyclic closure of $R$ as well as  the  multigraded deviations in terms of the
multigraded homotopy Lie algebra of $R$, and show that the LCM
lattice determines the Poincar\'e series of $R$. This is the
multigraded version of  Theorem 1 of
 \cite{Av02}.

Thanks are due to J.~Backelin and the faculty of the Department of
Mathematics in Stockholm, for useful discussions during my brief
stay in 1996.

\section{The  monomials of $Q_R(y,t)$ }\label{ymonomialsection}

The argument used in the proof of the next proposition has the
flavor of  \cite{Ba82} and uses Golod theory. It uses Lescot's
result, \cite{Le86}, that the Poincar\'e series $M$ over $R$
$$P_R^M(y,t) = \sum_{i=0}^{\infty}  \sum_{{\bf j}\ge 0}
y^{\bf j} {\dim}_k ({\Tor}_i^R(M,k)_{\bf j})
t^i={{f_R^M(y,t)}\over{Q_R(y,t)}}$$ where$f_r^M(y,t))\in
Z[y_1,\ldots, y_n][t]$.

\begin{proposition}\label{lcmcoefficients}
Let $Q_R(y,t)=1+\sum (\sum c_{\bj} y^{\bj})t^i$. Then $y^{\bj}$ is
equal to a least common multiple of a subset of the minimal
monomial generating set of $I$.
\end{proposition}
\begin{proof} By a standard argument, \cite{Ba82},  it is enough to prove  Proposition \ref{lcmcoefficients}
when the ideal $I$ is generated by squarefree monomials. In this
case the  $y$ monomials of $Q_R(y,t)$ are squarefree, \cite{Ba82}.
We do induction on the number of variables $n$ that divide the
monomials of  the minimal generating set of $I$, the case $n=1$
being trivial. Consider a term of $Q_R(y,t)$ whose coefficient in
$y$ does not involve $y_i$. We separate the generators of $I$ that
involve $x_i$. Let $I=L_i+ x_i J_i$ where  $L_i, J_i\subset
\Bbbk[x_1,\ldots,\hat{x_i}, \ldots, x_{n-1}]$ and let
$A=\Bbbk[x_1,\ldots, \hat{x_i}, \ldots, x_{n-1}]/L_i$ and
$A'=\Bbbk[x_1,\ldots, \hat{x_i}, \ldots, x_{n-1}]/(L_i+J_i)$.
Therefore $R\cong A[x_i]/x_i J_i$ and
 $A[x_i]\mapright{} R$ is Golod,   which implies that
$$P_R(y,t)={P_{A[x_i]}\over (1-t(P^R_{A[x_i]}-1))}=
{(1+ty_i)P_A\over (1-t^2P_{A[x_i]}^{x_iJ_i} ) }= {(1+ty_i)P_A
\over (1-t^2y_i P_A^{J_i})}=$$
$$={(1+ty_i)P_A \over (1-t y_i(P_A^{A'}-1))}=
{\prod (1+ty_i)\over Q_A+ty_iQ_A-ty_if_A^{A'}}.$$ The $y$
monomials that do not involve $y_i$ are terms of $Q_A$, so by
induction they have the desired form. Suppose now that a monomial
of $Q_R(y,t)$ involves all $y_i$: it is the product of all the
variables $y_i$ and is  equal to the least common multiple of all
the generators in $I$.
\end{proof}

In \cite{ChRe95} we noticed that the terms of $Q_R(y,t)$ when $I$
is a complete or almost complete intersection come from the set
$\{ (-1)^{l_\J}t^{(|\J|+l_\J)}m_\J\}$. Here $\J$ ranges over the
subsets of the minimal monomial generating set of $I$ and $l_\J$
is the number of connected components of the graph we get when we
connect the elements of $\J$ by an edge if they have a variable in
common. In \cite{ChRe95} we also confirmed that this is the case
when the Taylor resolution \cite{Ta60}, \cite{Ei97} of $S/I$ over
$S$   is minimal.

In Remark \ref{golod-denom-rem} we confirm that the same holds for
the monomial ideals $I$ such that $R/I$ is Golod. First we recall
and expand a comment by Fr\"oberg, \cite{Fr79}.

\begin{remark} Let $R=S/I$ be Golod and $K_\bullet$ be the Koszul
complex on the variables $\bar{x_i}$ of $R$. Then $H(K_\bullet)$
is isomorphic to $\Bbbk[X_{A_1},\ldots, X_{A_r}]/(X_{A_1},\ldots,
X_{A_r})^2$ where $A_i$ is a subset of the minimal monomial
generating set of $I$ and the degree of $X_{A_i}$ equals the
cardinality of $A_i$ .

Indeed since $K_\bullet$ admits trivial Massey operations the
product of two elements of $H_{\ge 1}(K_\bullet)$ is zero,
\cite{GuLe69}, \cite{Av98}. Moreover $H(K_\bullet)$ is always of
the form $\Bbbk[X_1, \ldots, X_r]/L$: $L$ is generated by
quadratics, the $X_i's$ form a minimal multiplicative generating
set for $H(K_\bullet)$, and $X_i X_j=(-1)^{|X_i| |X_j|} X_i X_j$
where $|X_i|$ is the degree of $X_i$. Let $F_\bullet$ be the
Taylor complex on the minimal monomial generators of $I$ over $S$,
\cite{Ta60} and let $E_\bullet$ be the minimal multigraded
resolution of $S/I$ over $Q$. $E_\bullet$ is a direct summand of
$F_\bullet$. It follows that the multidegrees of  the generators
of $E_i$ are among the multidegrees of the generators of $F_i$.
Therefore they are the multidegrees of  least common multiples of
subsets $A_j$ of $\I$. It follows that $H(K\bullet)=H(\Bbbk\otimes
E_\bullet)= \Bbbk[X_{A_1},\ldots, X_{A_r}]/L$. Since $X_{A_j}
X_{A_t}=0$,  $X_{A_j} X_{A_t} \in L$.
\end{remark}

\begin{remark}\label{golod-denom-rem}
Thus monomial Golod rings are special Golod rings. Let $R$ be a
Golod ring. From our previous remark and using the fact that
$m_A\neq m_{A_1}\cdots m_{A_t}$ for any partition of $A$ where
$X_A$ and $X_{A_j}$ are among the generators of $H(K_\bullet)$,we
see that $A$ is connected and $l_A=1$. Let $H_i(K)_{\bf j}$ denote
the subspace of $H_i(K)$ of multidegree $\bf j$. Thus the terms of
the denominator of the Poincar\'e series of $R$,
$Q_R(y,t)=1-\sum_i \sum_{\bf j} \dim_\Bbbk H_i(K)_{\bf j} y^{\bf
j} t^{i+1}$ of the Golod ring $R$ come from the set $\{
(-1)t^{(|\J|+1)}m_\J\}$ where $\J$ ranges over the subsets of the
minimal monomial generating set of $I$.
\end{remark}

Generic ideals and their minimal resolutions were first introduced
in \cite{BaPeSt98} and generalized in \cite{MiStYa00}. Let $I$ be
a generic ideal. If the multidegrees of  two minimal monomial
generators of $I$ are equal for some variable then there is a
third monomial generator of $I$ whose multidegree is strictly
smaller than the multidegree of the least common multiple of the
other two. The minimal resolution $E_\bullet$ of $S/I$ over $S$ is
determined by the Scarf complex of $\D_I$. The latter  consists of
the subsets $\J$ of the index set of $I$, $\Phi$, such that
$m_{\J}\neq m_{\L}$ for any other subset $\L$ of $\Phi$. We note
that $E_\bullet$ has the structure of a DG algebra,
\cite{BaPeSt98}. The Eagon resolution of $\Bbbk$ over $R=S/I$ was
discussed in \cite{GuLe69}. Below we give the generators for the
Eagon resolution when $I$ is generic.

\begin{theorem}
\label{res-thm} Let $R=S/I$ where $I$  is a generic ideal of $S$,
let $K_\bullet$ be the Koszul complex on the $x_i$ over $R$,
$X_0=0$, and  $X_i$ be the free $R$-module with generators $T_{\rm
L}$ where $|\L|=i$ and $\L$ is in $\Delta_I$,  the Scarf complex
of $I$. Then  $(Y_\bullet, d_\bullet)$ is a free resolution of
$\Bbbk$ over $R$, where
$$Y_i=\sum_{j+t_1+\cdots+t_l+l=i}
                  K_j \otimes
          X_{t_1}\otimes\cdots \otimes X_{t_l},$$
and
\begin{flushleft}
$d(T_\J)=(m_\J/x_\J ) e_\J,$
\end{flushleft}
\begin{flushleft}
$d(T_{\I_1}\otimes\cdots\otimes T_{\I_l})=
 d(T_{\I_1} \otimes \cdots\otimes T_{\I_{l-1}})
\otimes T_{\I_l}
  +$
\end{flushleft}
\begin{flushright}
$(-1)^{\displaystyle{\sum_{i=1}^{l-1}
 (|\I_i|+1)}}
{\displaystyle T_{\I_1}\otimes\cdots\otimes
  T_{\I_{l-2}}\otimes \delta(\I_{l-1},\I_l})(T_{\I_{l-1}}\wedge T_{\I_l}),$
\end{flushright}
where $\delta(\I_{l-1},\I_l)={{ m_{\I_{l-1}} m_{\I_l}\over
m_{\I_{l-1}\cup \I_l} }} $ if $\I_{l-1}\cup\I_l\in \Delta_I$ and 0
otherwise
\begin{flushleft}
$d(e_\L \otimes T_{\I_1}\otimes\cdots\otimes
T_{\I_l})d(e_\L)\otimes T_{\I_1}\otimes\cdots\otimes T_{\I_l}+
(-1)^{|\L|}e_\L\otimes d(T_{\I_1}\otimes\cdots\otimes
 T_{\I_l}).$
\end{flushleft}
\end{theorem}

\begin{proof} The theorem follows immediately from \cite{Iy97}, Theorem 1.2 and the preceeding comments.
\end{proof}

As an immediate corollary we get the characterization of Golod
rings when $I$ is a generic ideal, which generalizes Proposition
2.4 of \cite{Ga00}.

\begin{proposition} Let $I$ be a generic ideal. $S/I$ is Golod if and only if  $m_{\I_{1}} m_{\I_2}\neq m_{\I_{1}\cup \I_2} $ whenever  $\I_{1}\cup\I_2\in \Delta_I$.
\end{proposition}

\section{Multigraded deviations and acyclic closures}\label{multdeviations}

In this section we introduce the multigraded acyclic closure
$R\<X{\bigger}$ of $\Bbbk$ over $R$. For the construction of the
usual acyclic closure, we refer to Construction 6.3.1,
\cite{Av98}. To this end we choose $x_{\bj}$, the $\G$-variables
of $X_n$  so that $cls(\theta(x_{\bj}))$ has multidegree $\bj$ and
$\{ cls(\theta(x_{\bj}))\}$ is a  minimal multigraded generating
set  of $H_{n-1}(R\<X_{\le n-1}{\bigger})$.  We define the
multidegree of $x_{\bj}$ to be $\bj$ and we let $X_{n,\bj}$
consist of all $ x\in X_n$ whose  multidegree  is $ \bj$. We will
show that the cardinality of $X_{n,\bj}$ appears as an exponent in
a product decomposition of the multigraded Poincar\'e series.

First we state and prove the multigraded version of Remark 7.1.1,
\cite{Av98}.

\begin{proposition} For each formal multigraded power series with integer coefficients
$$P(y,t)=1+\sum_{i=1}^\infty
(\sum_{\bj} a_{i, \bj}y^{\bj}) t^i$$ where for each $i$,
$a_{i,\bj}=0$ for all but finitely many values of $\bj$,
 there exist uniquely defined integers
$e_{n,{\bj}}\in Z$ such that
$$P(y,t)= {{\prod_{i=1}^\infty \prod_{\bj}(1+y^{\bj} t^{2i-1})^{e_{2i-1,{ \bj}}} }\over
    {\prod_{i=1}^\infty \prod_{\bj}(1-y^{\bj} t^{2i})^{e_{2i,{\bj} }} }}
$$
and the product converges in the (t)-adic topology of the ring
$Z[y,t]$.
\end{proposition}

\begin{proof}
We set $P_0(t)=1$ and assume by induction that $P(y,t)\equiv
P_{n-1}(y,t) $ modulo $ t^{n+1}$. If $P(y,t)-P_{n-1}(y,t)\equiv
\sum e_{n,{\bj}}y^{\bj} t^n (\modulo\ t^{n+1})$, then  we set
$P_n(y,t)=P_{n-1}(y,t) \prod_{\bj} (1+y^{\bj} t^n)^{e_{n,\bj}}$ if
$n$ is odd and
 $P_n(y,t)=P_{n-1}(y,t) / \prod_{\bj} (1-y^{\bj} t^n)^{e_{n,\bj}}$ if $n$ is even.
Then it is clear that $P(t)\equiv P_n(t) (\modulo\ t^{n+1})$ and
the other assertions follow as well.
\end{proof}

 We define the {\it multigraded} $(n,{\bj})$ {\it deviation} of $R$, denoted by
$\e_{n,\bj}$, to be the exponent $e_{n, {\bj}}$ in the product
decomposition of the Poincar\' e series of $R$.  As in Theorem
7.1.3, \cite{Av98} it follows that
\begin{proposition} Let $I$ be a monomial ideal of $S$, $R=S/I$ and
$R\<X\bigger$ be the multigraded acyclic closure of $\Bbbk$ over
$R$. Then
$$P_R(y,t)= {{\prod_{i=1}^\infty \prod_{\bj}(1+y^{\bj} t^{2i-1})^{\card(X_{2i-1,{ \bj}})} }\over
    {\prod_{i=1}^\infty \prod_{\bj}(1-y^{\bj} t^{2i})^{\card(X_{2i,{\bj} }) }  }}
$$
and $\card X_{n,\bj}=\e_{n,\bj}(R)$.
\end{proposition}

Since $R\< X\bigger$ is multigraded, it follows that the homotopy
Lie algebra of $R$, $\pi^\bullet(R)=H Der_R^\gamma(R\< X{\bigger},
R\< X{\bigger})$ has a multigraded structure. Moreover $\pi^n(R)
\cong \Hom_{\Bbbk}({\Bbbk}X_n, \Bbbk)$, (Theorem 10.2.1,
\cite{Av98}), and $\pi^n_{\bj}(R) \cong
\Hom_\Bbbk({\Bbbk}X_{n,{-\bj}}, \Bbbk)$ while $\rank_\Bbbk
\pi^n_{\bj}(R)=\e_{n,-\bj}$.

\begin{remark}
We let $K_\bullet$ denote the Koszul complex on $x_1,\ldots, x_n$
over $S$ and $T_\bullet$ the Taylor complex on the minimal
generators of $I$ over $S$. The complexes $K_\bullet$, $T_\bullet$
have the structure of a DG $\Gamma$ algebra by Lemma 9,
\cite{Av02}; they are also multigraded  and
 $H_\bullet (T\otimes_S K)_{\bj}$
$\cong H_\bullet(T\otimes_S \Bbbk)_{\bj}\cong H_\bullet
(R\otimes_S K)_{\bj}$, where $\bj \in L_I$. In Lemma 11 of
\cite{Av02}, it is shown that $\pi^{\ge 2}(R)\cong \pi^\bullet
(R\otimes_S K)$.  We remark that $R\otimes_S K$ is a multigraded
DG $\Gamma$ algebra, and one can choose a multigraded DG $\Gamma$
algebra $U$ so that we have the factorization $R\otimes_S
K\mapright{} U\mapright{}\Bbbk$ with the properties of Lemma 11
and Remark 10 of \cite{Av02}: $\pi(R\otimes_S K)$ is the graded
$\Bbbk$ dual of  the residue of $H_\bullet(U\otimes_{R\otimes_S K}
\Bbbk)$ modulo multigraded relations and is multigraded. The
homomorphisms of Lemma 11, \cite{Av02} are multigraded and
$\pi^{l+2}_{\bj}(R)\cong \pi^{l}_{\bj} (R\otimes_S K)$. Finally we
have $\pi^{l}_{\bj} (T\otimes_S K)\cong \pi^{l}_{\bj} (T\otimes_S
\Bbbk)\cong \pi^{l}_{\bj}(R\otimes_S K$.
\end{remark}

Next we will make use of the LCM lattice $L_I$, \cite{GaPeWe99}
and the GCD graph $G_I$,\cite{Av02}. We recall that $L_I$ and
$G_I$ have vertices the least common multiples of the monomial
generators of $I$ while the edges of $G_I$ join least common
multiples that are relatively prime. In \cite{Av02} it is shown
that if $I$ and $I'$ are two monomial ideals of $S$ and $S'$
respectively with an isomorphism of lattices $\lambda:
L_I\mapright{} L_{I'}$ which induces an isomorphism of the GCD
graphs $\lambda: G_I\mapright{} G_{I'}$ then the homotopy Lie
algebras of $R=S/I$ and $R'=S/I'$ are isomorphic and the
Poincar\'e series of $R$ and $R'$ have the same denominator. We
will show that the multigraded version of this result is actually
true. We consider the multigraded acyclic closures $R<X>$ and
$R'<X'>$ of $\Bbbk$ over $R$ and $R'$.

\begin{remark}
Let $I\subset S=\Bbbk[x]$ and $I'\subset S'=\Bbbk[x']$ be two
monomial ideals, $K, K'$ the Koszul complexes on the $x_i$ and
$x_i'$ respectively and $\lambda: L_I\mapright{} L_{I'}$ an
isomorphism of lattices. If $\lambda$  induces an isomorphism of
the GCD  graphs $\lambda: G_I\mapright{} G_{I'}$ then $H(K\otimes
R)_{\bf j}\cong H(K'\otimes R')_{\lambda(\bf j)}$ through the
isomorphism that sends the algebra generators $X_A$ of $H(K\otimes
R)=\Bbbk[X_A]/L$ to the generators $X'_{\lambda(A)}$ of
$H(K'\otimes R')=\Bbbk[X'_{\lambda(A)}]/L'$, (here our notation is
as in the remarks preceeding Proposition 2.2). Indeed, let
$\lambda$ be as above. Since $\lambda$ is an isomorphism of the
LCM lattices the minimal resolution of $I$ gives rise to a minimal
resolution of $I'$ and the generators $X_A$ of  $H(K\otimes R)$
correspond to generators $X'_{\lambda(A)}$ of $H(K'\otimes R')$.
Since $\lambda$ is an isomorphism of the GCD graphs the relations
among  the $X_A$ correspond to relations among the
$X'_{\lambda(A)}$. Conversely if the Koszul homology algebras are
isomorphic as above then there is a an isomorphism of multigraded
vector spaces between the minimal resolutions of the two ideals
and a map $\lambda$ which is an isomorphism of lattices and GCD
graphs for the multidegrees that appear in the resolutions.
\end{remark}

\begin{proposition}
Let $I'$ be a monomial ideals of $S'=\Bbbk[x']$ such that
$\lambda: L_I\mapright{} L_{I'}$ is an isomorphism of lattices
that induces an isomorphism of the GCD  graphs $\lambda:
G_I\mapright{} G_{I'}$.  Then there is a map $\hat\lambda$ from
the set of the multidegrees of the variables $X_n$ to the set of
the multidegrees of the variables $X'_n$ such that $\card
X_{n,\bj}=\card X'_{n,\hat\lambda(\bj)}$ and $\e_{n,\bj}=\e'_{n,
\hat\lambda(\bj)}$.
\end{proposition}
\begin{proof} We define $\hat \lambda$ by induction on $n$.
Let $T'$ be the Taylor complex of $I'$. Then there is an
isomorphism of vector spaces $(T\otimes_S \Bbbk)_{\bj}\cong
(T'\otimes_{S'} \Bbbk)_{\lambda(\bj)}$, for  $\bj \in L_I$.
Suppose that $\hat \lambda$ is an extension of $\lambda$ and
$R[X_{\le n,\bj}]\otimes \Bbbk \cong R'[X'_{\le
n,\hat\lambda(\bj)}]\otimes \Bbbk$. Extend $\hat\lambda$ to the
multidegrees of the multigraded generators of $H_n(R[X_{\le n}])$
and $H_n(R'[X'_{\le n}])$. It follows that there is an isomorphism
of homotopy Lie algebras $\pi^{\bullet}_{-\bj} (T\otimes_S
\Bbbk)\cong \pi^{\bullet}_{-{\hat \lambda}(\bj)} (T'\otimes_{S'}
\Bbbk)$, and the proposition now follows.
\end{proof}

\begin{remark} Let $I$ be a monomial ideal and $I_{pol}$ be the square
free monomial ideal in $A=\Bbbk[z]$ that corresponds to $I$,
\cite{Fr82}. The map $\lambda $ that sends a monomial of $I$ to a
squarefree monomial in $I_{pol}$ is an isomorphism of lattices and
of GCD graphs. Moreover $Q_R(y,t)$ can be obtained from $Q_A(z,t)$
by applying $\lambda^{-1}$ to the monomial coefficients of $t$,
\cite{Ba82}.
\end{remark}

Theorem \ref{denomisom} now completes the multigraded version of
Theorem 1, \cite{Av02}.

\begin{theorem}\label{denomisom} Let $S=\Bbbk[x]$ and $S'=\Bbbk[x']$ be
polynomial rings over a field $\Bbbk$,  and $I\subset S$,
$I'\subset S'$ be ideals generated by monomials of degree at least
2. Let $\lambda: L_I\mapright{} L_{I'}$ and the induced $\lambda:
G_I\mapright{} G_{I'}$ be isomorphisms.  Suppose that  $Q_R(y,
t)=1+\sum (\sum c_{\bj} y^{\bj}) t^i$ then $Q_{R'}(y',t)= 1+\sum
(\sum c_{\bj} y^{\lambda(\bj)}) t^i$.
\end{theorem}

\begin{proof}
By Proposition \ref{lcmcoefficients} the multidegrees $\bj$ in
$Q_R(y,t)$ are in $L_I$, so $\lambda(\bj)$ makes sense. By the
above remark it is enough to examine the case where $I$ and $I'$
are both squarefree.  The key point is that since $I$ and $I'$ are
squarefree, $\lambda$ is  additive: $\lambda(\bf j+\bf
i)=\lambda(\bf j)+\lambda(\bf i)$. It follows that its extension
$\hat\lambda$ is also additive. Since
$$Q_R(y,t)={ \prod_{i=1}^n(1+ty_i)\over P_R(y,t)}=
 {{\prod_{i=1}^\infty \prod_{\bj}(1-y^{\bj} t^{2i})^{\rank\pi^{2i}_{- \bj}}  }
\over{\prod_{i=2}^\infty \prod_{\bj}(1+y^{\bj} t^{2i-1})^{\rank
\pi^{2i-1}_{-\bj}} }}
$$
the linear coefficient $c_{\bj}$ of the monomial $ y^{\bj}t^n$ of
$Q_R(y,t)$, depends on the deviations $\e_{i,\bff}$ when $i\le n$.
The corresponding combination of the deviations
$\e'_{i,\hat\lambda(\bff)}$ gives a term of $Q_{R'}(y',t)$ of
$t$-degree $n$. The theorem now follows.
\end{proof}

The condition that $\lambda$ is an isomorphism of GCD graphs is a
necessary condition as the example from \cite{Ga00} shows. The
isomorphism of LCM lattices for $I= (x_1^2, x_2^2x_3)$ and
$I'=(x_1x_2^2, x_1x_3^2)$ is not an isomorphism of GCD graphs and
$$Q_{S/I}(y,t)=1-t^2(y_1^2-y_2^2y_3)+t^4(y_1^2y_2^2y_3)$$
while
$$Q_{S/I'}(y,t)=1-t^2(y_1y_2^2+y_1y_3^2)-t^3(y_1y_2^2y_3^2).$$

\end{document}